\documentclass[10pt,oneside,a4paper]{jm}
\usepackage{amsfonts,amssymb,amsmath,stmaryrd} 

\newcommand{\abs}{\vskip 0.5em\noindent}
\newcommand{\Abs}{\paragraph{}\hspace{-1em}}
\newcommand{\AbsT}[1]{\paragraph{\hspace{-1em} #1}}

\newcommand{\Thm}{\AbsT{Theorem .}}
\newcommand{\ThmT}[1]{\AbsT{Theorem : #1.}}

\newcommand{\Conj}{\AbsT{Conjecture.}}

\newcommand{\bfa}{\noindent {\bf a)} }
\newcommand{\bfb}{\noindent {\bf b)} }

\newcommand{\N}{\mathbb N}

\newcommand{\F}{\mathbb F}

\newcommand{\cF}{\mathcal F}

\newcommand{\cK}{\mathcal K}
\newcommand{\cL}{\mathcal L}

\newcommand{\cM}{\mathcal M}
\newcommand{\cN}{\mathcal N}

\newcommand{\cO}{\mathcal O}

\newcommand{\cW}{\mathcal W}

\newcommand{\fA}{\mathfrak A}

\newcommand{\al}{\alpha}

\newcommand{\bt}{\beta}
\newcommand{\gm}{\gamma}

\newcommand{\dt}{\delta}
\newcommand{\Dt}{\Delta}
\newcommand{\eps}{\epsilon}

\newcommand{\lb}{\lambda}

\newcommand{\om}{\omega}

\newcommand{\Om}{\Omega}

\newcommand{\Aut}{\text{Aut}}

\newcommand{\GL}{\text{GL}}

\newcommand{\Hom}{\text{Hom}}

\newcommand{\Ind}{\text{Ind}}

\newcommand{\Inn}{\text{Inn}}

\newcommand{\Br}{\text{Br}}
\newcommand{\modcat}[2]{\textbf{mod}_{#1}\text{-}#2}

\newcommand{\Out}{\text{Out}}

\newcommand{\Res}{\text{Res}}

\newcommand{\Stab}{\text{Stab}}

\newcommand{\Tr}{\text{Tr}}

\newcommand{\ld}{,\ldots\hskip0em ,}
\newcommand{\lr}[1]{\langle #1\rangle}

\newcommand{\wti}[1]{\widetilde{#1}}

\newcommand{\mt}{\mapsto}

\newcommand{\ra}{\rightarrow}

\newcommand{\Lra}{\leftrightarrow}
\newcommand{\lra}{\longrightarrow}
\newcommand{\lla}{\longleftarrow}

\newcommand{\cn}{\colon}

\newcommand{\spmid}{\,\mid\,}

\newcommand{\sseq}{\subseteq}

\newcommand{\emp}{\emptyset}

\newcommand{\bsl}{\backslash}
\newcommand{\tm}{\times}
\newcommand{\otm}{\otimes}

\newcommand{\GAP}{{\sf GAP}}
\newcommand{\MA}{{\sf MeatAxe}}

\begin{document}
\raggedbottom
\pagestyle{myheadings}
\thispagestyle{empty}
\begin{center} 
{\Large\bf
Brou{\'e}'s abelian defect group conjecture \vspace*{0.2em} \\
for the sporadic simple Janko group $J_4$ revisited  
\vspace*{1em}
}
        {\bf\large
        {Shigeo Koshitani}$^{\text a,*}$, 
        {J{\"u}rgen M{\"u}ller}$^{\text b}$,
        {Felix Noeske}$^{\text c}$  
               } 
\footnote{
$^*$ Corresponding author. \\
\indent {\it E-mail addresses:} 
koshitan@math.s.chiba-u.ac.jp (S.~Koshitani), \\
juergen.mueller@math.rwth-aachen.de (J.~M{\"u}ller), \\
Felix.Noeske@math.rwth-aachen.de (F.~Noeske).  
\\
{\bf \today}
}

\medskip

%
{\it
${^{\mathrm{a}}}$Department of Mathematics, Graduate School of Science, \\
Chiba University, Chiba, 263-8522, Japan \\
${^{\mathrm{b, \, c}}}$Lehrstuhl D f{\"u}r Mathematik,
RWTH Aachen University, 52062 Aachen, Germany}
\end{center}

\begin{center}
{\sf Dedicated to Geoffrey R. Robinson on the occation of 
his sixtieth birthday}
\end{center}

%
%

\hrule 
\begin{abstract} 

\noindent
We show that the $3$-block of the sporadic simple Janko group $J_4$
with defect group $C_3\tm C_3$, and the principal $3$-block of the 
alternating group $\fA_8$ are Puig equivalent, answering a question
posed in \cite{KosKunWak}. To accomplish this, we apply computational
techniques, in particular an explicit version of the Brauer construction.
\end{abstract}
\hrule 

\section{Introduction}\label{intro}

\abs
In recent years, much impetus in modular representation theory 
of finite groups has originated from attempts to prove
various fascinating deep conjectures. Two of them 
are Brou{\'e}'s Abelian Defect Group Conjecture \cite{Bro1990}
and a strengthening, Rickard's Splendidness Conjecture \cite{Ric},
which for the purpose of the present paper may be stated as follows:

\Conj
Let $k$ be an algebraically closed field of characteristic $p>0$,
and let $G$ be a finite group. Let $A$ be a block of $kG$
having an \emph{abelian} defect group $P$, let $N_G(P)$ be the 
normaliser of $P$ in $G$, and let $B$ be the
block of $k[N_G(P)]$ which is the Brauer correspondent of $A$.
Then Brou{\'e}'s Conjecture says that $A$ and $B$ are derived
equivalent, and Rickard's Conjecture says that there even 
is a splendid derived equivalence between $A$ and $B$.

\abs\abs
In general, Brou{\'e}'s and Rickard's Conjectures currently
are widely open. They have been proven for a number of cases, 
where for an overview we refer to \cite{ChuRic}. In particular,
in \cite{KosKun} it is shown that both conjectures hold true
whenever $A$ is a principal block having a defect group
$P\cong C_3\tm C_3$ isomorphic to the elementary abelian group 
of order $9$. This moves non-principal blocks with defect group
$P$ into the focus of interest, where, 
in view of the successful reduction strategy for principal blocks
used in \cite{KosKun}, and a possible, as yet non-existent,  
generalisation to non-principal blocks, it seems worthwhile 
to proceed with blocks of quasi-simple groups.
There are a few results already known, see for example
\cite{KosKunWak2002, KosKunWak2004, KosKunWak,
KosMue, KosMueNoe2012, MueSch},
which indicate that fairly often a non-principal 
block with defect group $P$ is actually Morita 
equivalent to a principal block of a different (smaller) group.
The present paper is another step in clarifying this relationship:

\abs
Letting $p:=3$, here we consider the sporadic simple Janko group $J_4$.
Then $k[J_4]$ has a unique block $A$ of defect $2$, hence having a defect
group $P\cong C_3\tm C_3$.
In order to verify Brou{\'e}'s and Rickard's Conjectures for $A$,
it is shown in \cite{KosKunWak} that $A$ is actually
Morita equivalent to the principal block $A'$ of $k[\fA_8]$, 
also having defect group $P$,
where $\fA_8$ is the alternating group on $8$ letters.
But it is left open, see \cite[Question 6.14]{KosKunWak},
whether or not $A$ and $A'$ are even Puig equivalent.
We are now able to answer this affirmatively:

\Thm\label{puigequiv}
The blocks $A$ and $A'$ are Puig equivalent.

\abs\abs
This also sheds some light on Rickard's Splendidness Conjecture 
for $A$ and its Brauer correspondent $B$ in $N_{J_4}(P)$, where in 
\cite[(6.13)]{KosKunWak} an indirect proof is given, running as follows:
Letting $B'$ be the Brauer correspondent of $A'$ in $N_{\fA_8}(P)$,
since both $A$ and $A'$, as well as $B$ and $B'$ are mutually Morita 
equivalent, and Green correspondences between $A'$ and $B'$, as well
as $A$ and $B$ are suitably related, the proof of Brou{\'e}'s Conjecture 
for $A'$ given in \cite[Example 4.3]{Oku}, which by \cite[Corollary 2]{Oku2000} 
also proves Rickard's Conjecture, works entirely similar for $A$.
Now, our approach provides explicit Puig equivalences between both 
$A$ and $A'$, as well as $B$ and $B'$, hence these directly transport a 
splendid tilting complex between $A'$ and $B'$ to such a complex 
between $A$ and $B$, thus Rickard's Conjecture for $A$ is completely
reduced to $A'$.

\AbsT{Strategy.}
Our strategy 
is to rework the approach in \cite{KosKunWak}
explicitly. This means, we fix a concrete realisation of $J_4$,
suitable to work with computationally, and a certain configuration 
of subgroups, in particular containing a copy of $\fA_8$;
this is carried out in Section \ref{setting}.
While doing so, we painstakingly take care that all choices
made are unique up to simultaneous $J_4$-conjugacy.

\abs
Having this in place, since the functors provided in \cite{KosKunWak}
inducing stable equivalences between $A$ and $B$, as well as
$A'$ and $B'$ coincide with Green correspondence on simple
modules, they can, in principle at least, be evaluated explicitly
on simple modules. But the Morita equivalence between $B$ and $B'$
used in \cite{KosKunWak} is based on the abstract theory of blocks
with normal defect groups, see \cite[Theorem A]{Kuel} and \cite[Proposition 14.6]{Pui},
hence is replaced here by explicit functors
relying on our fixed configuration of subgroups of $G$;
this is carried out in Section \ref{stroke}.
We would like to mention that we have been led to consider the block $A$
while preparing our earlier paper \cite{KosMueNoe2012}, 
on the double cover of the sporadic simple Higman-Sims group 
(which is much smaller than $J_4$), where a similar, but subtly 
different problem occurs in the analysis of local subgroups.

\abs
Alone, simple $A$-modules are much too large to be dealt
with directly by an explicit approach. (This makes up 
a decisive difference to our earlier considerations in \cite{KosMueNoe2012}.)
Instead, we use the Brauer construction, applied to $p$-permutation
modules, to facilitate the explicit determination of
Green correspondents of certain simple $A$-modules.
The underlying theory is presented in Section \ref{brauer},
which is largely based on \cite{Bro}, with a view towards
explicit computations.

\abs
A few comments on this computational approach seem to be in order:
In practice, using this technique we are able to reduce the size 
of objects to be handled computationally dramatically, as is seen
for example in our main application in \ref{applybrauer}.
From a more conceptual point of view, the description in 
\ref{brauertrivsrc} shows a formal similarity to so-called 
\emph{condensation} of permutation modules, a well-known workhorse
in computational representation theory, see for example 
\cite[Sections 9, 10]{MueHabil}.
Condensation is formally described as an application of a 
Schur functor associated with a $p'$-subgroup $K$.
(The use of the letter `$K$' in \ref{brauerperm} is reminiscent of 
the German writing of `Kondensation'.) Here, the role of $K$ is eventually 
played by the $p$-subgroup $P$, so that in a sense we are dealing 
with a `$p$-singular' generalisation of condensation.
Moreover, the Brauer construction features prominently in modular 
representation theory of finite groups; we only mention the 
comments in \cite[Section 4]{Ric} concerning its application to
splendid tilting complexes. Thus, due to its general nature
and the gain in computational efficiency,
we are sure that this technique will face more applications.

\Abs
In order to facilitate the necessary computations, 
we make use of the computer algebra system \GAP{} \cite{GAP}, 
to deal with finite groups, in particular
permutation and matrix groups, and with ordinary and Brauer
characters of finite groups. In particular, we make use of the
character table library \cite{CTblLib} of \GAP{},
the \GAP-interface \cite{AtlasRep} to the database \cite{ModAtlasRep},
and the SmallGroups library \cite{SmGrp} of \GAP{}.
Moreover, we use the  computer algebra system \MA{} \cite{Par,MA}, 
and its extensions 
\cite{LuxMueRin,  LuxSzoke, LuxSzokeII, LuxWie} to deal with
matrix representations over finite fields.

\abs
We remark that, although for the theoretical developments
we fix an algebraically closed field $k$ of characteristic $3$,
explicit computations can always be done of a suitably
large finite field of characteristic $3$, where for
the computations to be described here even the prime field $\F_3$
turns out to be large enough.
We assume the reader familiar with the relevant concepts
of modular representation theory of finite groups and of
the various notions of equivalences between categories, 
as general references see \cite{KoeZim,NT,The};
our group theoretical notation is standard, borrowed from \cite{Atlas}.

\Abs{\bf Notation.}
Throughout this paper a module means a finitely generated 
{\it right} module unless stated otherwise.
Let $G$ be a finite group. 
We denote by $1_G$ the trivial character of $G$, and
by $k_G$ the trivial $kG$-module.
We write $H \leq G$ 
when $H$ is a subgroup of $G$,
and $H < G$ when $H$ is a proper subgroup of $G$.
Let $H \leq G$, and let $V$ and $W$ be a $kG$-module and 
a $kH$-module, respectively.
Then we write $\Res^G_H(V)$ or $\Res_H(V)$ 
for the restriction of $V$ to $H$, and
$\Ind_H^G(W)$ or $\Ind^G(W)$ 
for the induced module (induction) $W \otimes_{kH}kG$ of $W$ to $G$.
For a block algebra $B$ of $kG$, 
we denote by $1_B$ the block idempotent of $B$.
For a $kG$-module $X$ we write $X^\vee$ for the $k$-dual of $X$, namely
$X^\vee := \Hom_k(X, k)$.
For a subset $S$ of $G$ and an element $g \in G$ we write $S^g$ for $g^{-1}Sg$.
We write $\mathbb N$ and $\mathbb N_0$ for the sets of all positive integers
and of all non-negative integers, respectively.
Let $\mathfrak S_n$ and $\mathfrak A_n$ denote
the symmetric group and the alternating group on $n$ letters, respectively.


\section{Brauer construction}\label{brauer}

\AbsT{Brauer construction.}
Let $k$ be a field of characteristic $p>0$, and let $G$ be a finite group. 
For any finitely generated (right) $kG$-module $V$ and 
any subgroup $K\leq G$ let $V^K:=\{v\in V;vg=v\text{ for all }g\in K\}$ 
be the set of $K$-fixed points in $V$. Hence
$V^K$ naturally becomes a $k[N_G(K)]$-module, which since 
$K$ acts trivially induces the structure of a $k[N_G(K)/K]$-module on $V^K$.
Moreover, the $k$-linear trace map on $V$ is defined as 
$\Tr_K^G\cn V^K\ra V^G\cn v\mt\sum_{g\in K\bsl G}vg$,
the sum running over a system of 
representatives $g$ of the right cosets $Kg$ of $K$ in $G$.

\abs
Now, if $P\leq G$ is a $p$-subgroup, then let
$$ V(P):=V^P/\left(\sum_{Q<P}\Tr_Q^P(V^Q)\right) ,$$
where the sum runs over all proper subgroups $Q$ of $P$.
Since for all $g\in N_G(P)$ and $Q\leq P$ we have
$(\Tr_Q^P(v))^g=(\Tr_{Q^g}^P(v^g))\in V^P$, for all $v\in V^Q$,
we conclude that $\sum_{Q<P}\Tr_Q^P(V^Q)\leq V^P$ is a
$k[N_G(P)]$-submodule, hence $V(P)$ becomes a $k[N_G(P)]$-module as well.
Moreover, the natural map $\Br^P\cn V^P\ra V(P)$, 
being called the associated Brauer map, is a homomorphism of 
$k[N_G(P)]$-modules; note that in particular $\Br^P$ commutes
with taking direct sums.

\AbsT{$p$-permutation modules.}\label{brauertrivsrc}
Recall that $V$ is called a $p$-permutation $kG$-module
if $\Res_P(V)$ is a permutation $kP$-module for any $p$-subgroup $P\leq G$,
and that this equivalent to saying that all indecomposable direct summands
of $V$ are trivial source $kG$-modules, see \cite[Proposition 27.3]{The}.

\abs\bfa
Here, we are content with much less, namely we consider the case 
where $\Res_P(V)$ is a permutation $kP$-module for some fixed 
$p$-subgroup $P\leq G$.
The following facts are well-known, and stated for example
in \cite[Section 1]{Bro} and \cite[Exercise 11.4]{The}. Still, we
include the details which will be needed explicitly later on:

\abs
Let $\Om$ be a $P$-stable $k$-basis of $V$; we also write $V=k[\Om]$.
Let $\Om=\coprod_{i=1}^r\Om_i$ be the decomposition of $\Om$ into
$P$-orbits, for some $r\in\N_0$, where we assume that
$$ \Om^P:=\{\om\in\Om;\om g=\om\text{ for all }g\in P\}
=\coprod_{i=1}^s\Om_i=\coprod_{i=1}^s\{\om_i\}, $$ 
%
%
for some $s\in\{0 \ld r\}$, is the set of $P$-fixed points in $\Om$.
Then we have $\Res_P(V)=\bigoplus_{i=1}^r k[\Om_i]$ as $kP$-modules,
and $k[\Om_i]^P=k[\Om_i^+]$ where $\Om_i^+:=\sum_{\om\in\Om_i}\om$
is the associated orbit sum, for all $i\in\{1\ld r\}$.
Thus $\{\Om_1^+\ld\Om_r^+\}$ is a $k$-basis of $V^P$,
and as $k$-vector spaces we get
$$ V(P)=\Br^P(V^P)\cong\bigoplus_{i=1}^r\Br^P(k[\Om_i]^P) 
  =\bigoplus_{i=1}^r \Br^P(k[\Om_i^+]) .$$
If $\om_0\in\Om_i$ for some $i>s$, then for 
$Q:=\Stab_P(\om_0)<P$ we get 
$\Tr_Q^P(\om_0)=\sum_{\om\in\Om_i}\om=\Om_i^+$,
implying that $\Br^P(k[\Om_i^+])=0$; for $i\leq s$ we have $\Om_i=\{\om_i\}$
and thus $\Tr_Q^P(\Om_i^+)=\Tr_Q^P(\om_i)=[P\cn Q]\cdot\om_i=0$,
for all $Q<P$, showing that $\Br^P$ is the identity map on 
$k[\Om_i^+]=k[\om_i]$, hence so is on $k[\Om^P]$.
Thus we conclude that the Brauer map induces a $k$-vector space
isomorphism $\Br^P\cn k[\Om^P]\ra V(P)$.

\abs\bfb
Now we additionally assume that $\Om$ is $N_G(P)$-stable, 
that is $\Res_{N_G(P)}(V)$ is a permutation $k[N_G(P)]$-module 
with respect to $\Om$; note that this in particular
holds whenever $V$ is a permutation $kG$-module.
Then it follows that $N_G(P)$ permutes
the $P$-orbits in $\Om$ of any fixed length amongst themselves.
Hence both $\Om^P$ and 
$\Om - \Om^P$ are $N_G(P)$-stable.
Since we have already seen that 
%
%
$k[\Om - \Om^P]=\sum_{Q<P}\Tr_Q^P(V^Q)$, as $k[N_G(P)]$-modules we get
$$ V^P=k[\Om^P]\oplus k[\Om - \Om^P]
      =k[\Om^P]\oplus\left(\sum_{Q<P}\Tr_Q^P(V^Q)\right) .$$
Hence we conclude that the the induced map 
$\Br^P\cn k[\Om^P]\ra V(P)$ is even an isomorphism of $k[N_G(P)]$-modules;
in particular, $V(P)$ is a permutation $k[N_G(P)]$-module. Note that
this implies that $V(P)$ is a $p$-permutation $k[N_G(P)]$-module 
whenever $V$ is a $p$-permutation $kG$-module, see \cite[(3.1)]{Bro}.

\AbsT{Permutation modules.}\label{brauerperm}
Let $H\leq G$ be a subgroup, and let $k[\Om]$ be the
permutation $kG$-module associated with the set
$\Om:=H\bsl G$ of right cosets of $H$ in $G$; 
we also write $\om_g\in\Om$ for the coset $Hg\sseq G$, where $g\in G$.
Moreover, let $K\leq G$ be an arbitrary subgroup. We proceed to derive
a description of $\Om^K$, and of the structure of $k[\Om^K]$
as a $k[N_G(K)]$-module.
The following results are also shown in \cite[(1.3), (1.4)]{Bro}, 
in the case of $K=P$ being a $p$-group, using Mackey's Formula and 
Higman's Criterion. Although only the latter case will be
relevant in our applications, we present a general straightforward 
proof, in the spirit of the explicit approach taken here:

\abs
Since for any $g\in G$ we have $\Stab_G(\om_g)=H^g$, we conclude that
$\Om^K\neq\emp$ if and only if $K\leq H^g$ for some $g\in G$.
Hence, if this is not the case then we have $(k[\Om])(K)=\{0\}$ anyway,
thus, by going over to some $G$-conjugate of $H$ if necessary,
we may assume that $K\leq H$;
note that replacing $H$ like this just amounts to going over
to an equivalent permutation action of $G$.

\abs
Let $K_1\ld K_t\leq H$, for some $t\in\N$, be a set of 
representatives of the $H$-conjugacy classes of subgroups 
of $H$ being $G$-conjugate to $K$, and let $g_i\in G$ 
such that $K_i^{g_i}=K$, for $i\in\{1\ld t\}$, where 
we may assume that $K_1=K$ and $g_1=1$. 
Now, for $g\in G$ we have $\om_g\in\Om^K$ if and only if $K^{g^{-1}}\leq H$.
If this is the case, then there are $i\in\{1\ld t\}$ and $h\in H$ such that
$K^{g^{-1}}=K_i^{h^{-1}}=K^{g_i^{-1}h^{-1}}$,
implying $g_i^{-1}h^{-1}g\in N_G(K)$, and thus $g\in Hg_iN_G(K)\sseq G$.
Conversely, if $g=hg_i n\in Hg_iN_G(K)$, for some $h\in H$ and 
$n\in N_G(K)$, then we have 
$K^{g^{-1}}=K^{n^{-1}g_i^{-1}h^{-1}}=K_i^{h^{-1}}\leq H$.
Hence in conclusion we get
$$ \Om^K=\{\om_{g_in}\in\Om;i\in\{1\ld t\},n\in N_G(K)\} .$$

\abs
Next we note that the double cosets 
$Hg_1N_G(K)\ld Hg_rN_G(K)\sseq G$ are pairwise distinct:
Assume that there are $i\neq j\in\{1\ld t\}$ such that
$hg_in=h'g_jn'$, for some $h,h'\in H$ and $n,n'\in N_G(K)$; then 
$$ K_i^{h^{-1}h'}=K_i^{g_inn^{\prime -1}g_j^{-1}}
  =K^{nn^{\prime -1}g_j^{-1}}=K^{g_j^{-1}}=K_j $$
shows that $K_i$ and $K_j$ are $H$-conjugate, a contradiction.
Thus the above description of $\Om^K$ shows that 
we have $k[\Om^K]=\bigoplus_{i=1}^t k[\om_{g_i}N_G(K)]$ 
as permutation $k[N_G(K)]$-modules.
Moreover, from
$$ \Stab_{N_G(K)}(\om_{g_i})=N_G(K)\cap H^{g_i}=N_{H^{g_i}}(K)
   =N_{H^{g_i}}(K_i^{g_i})=N_H(K_i)^{g_i} $$
and $N_G(K)=N_G(K)^{g_i}$, we get
$$ k[\Om^K]\cong\bigoplus_{i=1}^t k[N_{H^{g_i}}(K)\bsl N_G(K)] 
           \cong\bigoplus_{i=1}^t k[N_H(K_i)\bsl N_G(K_i)]^{g_i} ,$$
where the latter summands denote the $k[N_G(K)]$-modules
obtained from the $k[N_G(K_i)]$-modules $k[N_H(K_i)\bsl N_G(K_i)]$
by transport via $g_i$.

\abs
This completes our description of $k[\Om^K]$ as $k[N_G(K]$-module.
Note that 
the 
$i$-th summand in the above decomposition is 
the trivial $k[N_G(K)]$-module if and only if we have 
$N_G(K_i)=N_H(K_i)$, or equivalently $N_G(K_i)\leq H$. 
Finally, we remark that going over to cardinalities in particular yields
$|\Om^K|=\sum_{i=1}^t\frac{|N_G(K)|}{|N_H(K_i)|}$,
showing that the above result is a special case 
of the `induction formula for marks' (in the theory
of Burnside rings) given in \cite[Theorem 2.2]{Pfe}.

\abs
Underlying our application of the Brauer construction 
now is the following

\ThmT{Brou\'e-Puig \cite[(3.4)]{Bro}}\label{puigthm}
Let $V$ be an indecomposable trivial source $kG$-module having
$P$ as a vertex, and let $f$ be the Green correspondence 
with respect to the triple $(G,P,N_G(P))$, see \cite[Theorem 4.4.3]{NT}. 
Then the $k[N_G(P)]$-module $V(P)$ coincides with the 
Green correspondent $f(V)$ of $V$.


\section{The setting}\label{setting}

\abs
From now on let $k$ be an algebraically closed field of characteristic $3$.
Moreover, we fix a $3$-modular system $(\cK,\cO,k)$ which is large enough.
That is to say, 
$\mathcal O$ is a complete discrete valuation ring of
rank one such that its quotient field $\mathcal K$ is
of characteristic zero, and its residue field
$k=\mathcal O/\mathrm{rad}(\mathcal O)$ is of
characteristic $3$, and that $\mathcal K$ and $k$ are
splitting fields for all
the (finitely many) groups occurring in the sequel.

\AbsT{The group $J_4$.}\label{j4}
Let $G:=J_4$ be the sporadic simple Janko group $J_4$. 
Then using the ordinary character table of $G$ given
in \cite[p.190]{Atlas}, also available electronically 
in the character table library of \GAP{},
it follows that $kG$ has a unique block $A$ having a
defect group $P$ of order $9$; hence $P\cong C_3\tm C_3$.

\abs
Moreover, $G$ has a unique conjugacy class consisting
of elements of order $3$, and if $Q<P$ is a subgroup of
order $3$ then we have $N_G(Q)\cong 6.M_{22}.2$, see \cite[Section 3]{KleWil}. 
Now $6.M_{22}.2$ has precisely two conjugacy classes consisting
of elements of order $3$, with associated centralisers of shape
$6.M_{22}$ and $(2^3\tm 3^2).2$, respectively, see \cite[p.41]{Atlas}.
Thus we infer that $G$ has a unique conjugacy class of subgroups 
isomorphic to $C_3\tm C_3$, and we have $C:=C_G(P)\cong 2^3\tm P$.

\abs
By \cite[Theorem 4.5]{Wak} the decomposition matrix of $A$ is as given in 
Table \ref{j4dectbl}, where the irreducible ordinary characters
belonging to $A$ are numbered as in \cite[p.190]{Atlas},
and their degrees are recorded as well. The simple $A$-modules
are just named $S_1\ld S_5$; their dimensions are immediately
read off the decomposition matrix.

\begin{table}\caption{The decomposition matrix of $A$.}\label{j4dectbl}
$$ \begin{array}{|r|l|ccccc|}
\hline
\text{degree} & A & S_1 & S_2 & S_3 & S_4 & S_5 \\
\hline
\hline
     4\,290\,927 & \chi_{14} & 1 & . & . & . & . \\
    95\,288\,172 & \chi_{21} & . & 1 & . & . & . \\
   300\,364\,890 & \chi_{25} & . & . & 1 & . & . \\ 
   393\,877\,506 & \chi_{27} & . & . & . & 1 & . \\ 
   394\,765\,284 & \chi_{28} & 1 & . & . & . & 1 \\ 
   493\,456\,605 & \chi_{30} & 1 & 1 & . & 1 & . \\ 
   690\,839\,247 & \chi_{31} & . & . & 1 & . & 1 \\ 
   789\,530\,568 & \chi_{35} & . & 1 & 1 & 1 & . \\ 
1\,089\,007\,680 & \chi_{41} & 1 & . & 1 & 1 & 1 \\ 
\hline
\end{array} $$

\end{table}

\AbsT{The group $\fA_8$.}\label{a8}
Let $G':=\fA_8$ be the alternating group on $8$ letters,
and let $P\cong C_3\tm C_3$ be a Sylow $3$-subgroup of $G'$;
we will see in \ref{subgrps} that we may safely reuse the letter `$P$' here.
Hence we have $C_{G'}(P)=P$, and $G'$ has precisely 
two conjugacy classes consisting of elements of order $3$, 
of cycle type $[3^2,1^2]$ and $[3,1^5]$, respectively, 
where $P$ contains four elements from each of these classes. 
Hence $G'$ has precisely two conjugacy classes of subgroups 
of order $3$. 

\abs
Let $A'$ be the principal block of $kG'$. 
Then the decomposition matrix of $A'$ is as given in Table \ref{a8dectbl},
where the simple $A'$-modules $1,7,13,28,35$ are denoted 
by their dimensions, and similarly the irreducible 
ordinary characters $\chi_1\ld\chi_{70}$ belonging to $A'$ 
are indexed by their degrees. The irreducible ordinary 
and Brauer characters belonging to $A'$ can be found in
\cite[p.22]{Atlas} and \cite[p.49]{ModularAtlas}, respectively,
and are available electronically via \cite{ModularAtlasProject}
as well as in the character table library of \GAP{}.

\begin{table}\caption{The decomposition matrices of $A'$ and $B'$.}
\label{a8dectbl}
$$ \begin{array}{|l|ccccc|}
\hline
A' & 1 & 7 & 13 & 28 & 35 \\
\hline
\hline
\chi_{1}  & 1 & . & . & . & . \\
\chi_{7}  & . & 1 & . & . & . \\
\chi_{14} & 1 & . & 1 & . & . \\ 
\chi_{20} & . & 1 & 1 & . & . \\ 
\chi_{28} & . & . & . & 1 & . \\ 
\chi_{35} & . & . & . & . & 1 \\ 
\chi_{56} & 1 & 1 & 1 & . & 1 \\ 
\chi_{64} & 1 & . & . & 1 & 1 \\ 
\chi_{70} & . & 1 & . & 1 & 1 \\ 
\hline
\end{array}
\quad\quad\quad\quad\quad\quad\quad\quad\quad\quad\quad
\begin{array}{|l|ccccc|}
\hline
B' & 1a & 1b & 1c & 1d & 2 \\
\hline
\hline
\chi_{1a} & 1 & . & . & . & . \\
\chi_{1b} & . & 1 & . & . & . \\
\chi_{1c} & . & . & 1 & . & . \\ 
\chi_{1d} & . & . & . & 1 & . \\ 
\chi_{2}  & . & . & . & . & 1 \\ 
\chi_{4a} & 1 & . & 1 & . & 1 \\ 
\chi_{4b} & . & 1 & . & 1 & 1 \\ 
\chi_{4c} & 1 & . & . & 1 & 1 \\ 
\chi_{4d} & . & 1 & 1 & . & 1 \\ 
\hline
\end{array} $$ 
\end{table}

\AbsT{Local subgroups of $G'$.}\label{hprime}
Let $H':=N_{G'}(P)\cong P\cn D_8$, see \cite[p.22]{Atlas}. 
Since $P$ is abelian, implying that $H'$ controls $G'$-fusion in $P$,
we deduce that $D_8$ permutes transitively the $4$-sets 
of non-trivial elements of $P$ of fixed cycle type, hence $H'$ has
precisely two conjugacy classes of subgroups of order $3$. 

\abs
Using the facilities of \GAP{} to deal with groups and their characters, 
we determine the ordinary character table of $H'$. It is given in Table 
\ref{hprimectbl}, where again the irreducible characters are indexed 
by their degrees, conjugacy classes are denoted by the orders of the elements 
they contain, and centraliser orders are recorded as well.
We choose notation such that $\chi_{1a}$ is the trivial character, 
and that $\chi_{1b}$ is distinguished amongst the three non-trivial 
linear characters, for example by having an element of order $4$ 
in its kernel.
(Of course, the characters of $H'\cong P\cn D_8$ can also be
determined easily from those of $P$ and $D_8$ via Clifford theory, 
see \cite[Chapter 3.3]{NT}. But since we need the character table 
explicitly anyway, we found a direct computation appropriate here.)

\abs
Now a computation with \GAP{} shows that
$\Inn(H')\cong H'$ is a normal subgroup of index $2$ in $\Aut(H')$,
and that any non-inner automorphism $\om\in\Aut(H')$ induces a non-inner 
automorphism of $H'/P\cong D_8$ and interchanges the two 
$H'$-conjugacy classes of subgroups of order $3$ of $P$.
In particular, $\om$ induces a table automorphism of the
ordinary character table of $H'$, fixing $\chi_{1a}$, $\chi_{1b}$, 
and $\chi_2$, but interchanging 
$\chi_{1c}\Lra\chi_{1d}$, $\chi_{4a}\Lra\chi_{4c}$, and 
$\chi_{4b}\Lra\chi_{4d}$.

\abs
Let $B'$ be the principal block of $kH'$, that is the
Brauer correspondent of $A'$.
Since $A'$ is the only block of $kG'$ having maximal
defect, by Brauer's First Main Theorem , see \cite[Theorem 5.2.15]{NT},
we conclude that $B'$ is the only block of $kH'$, that is $B'=kH'$. 
Since the irreducible Brauer characters
of $H'$ are in bijection with the irreducible ordinary characters
of $H'/P\cong D_8$, the decomposition matrix of $B'$ is immediate
from the ordinary character table of $H'$. It is given in
Table \ref{a8dectbl}, where again the simple $B'$-modules
are denoted by their dimensions. Note that the non-inner automorphism
$\om$ fixes $1a$, $1b$, and $2$, but interchanges $1c\Lra 1d$. 

\begin{table}\caption{The ordinary character table of $H'$.}\label{hprimectbl}
$$ \begin{array}{|l|rrrrrrrrr|}
\hline
\text{centraliser} & 72&  8& 12& 12& 18& 18&  4&  6&  6 \\
\text{class}       & 1A& 2A& 2B& 2C& 3A& 3B& 4A& 6A& 6B \\
\hline
\hline
\chi_{1a} & 1& 1& 1& 1& 1& 1& 1& 1& 1 \\
\chi_{1b} & 1& 1&-1&-1& 1& 1& 1&-1&-1 \\
\chi_{1c} & 1& 1&-1& 1& 1& 1&-1&-1& 1 \\
\chi_{1d} & 1& 1& 1&-1& 1& 1&-1& 1&-1 \\
\chi_2    & 2&-2& 0& 0& 2& 2& 0& 0& 0 \\
\chi_{4a} & 4& 0& 0& 2&-2& 1& 0& 0&-1 \\
\chi_{4b} & 4& 0& 0&-2&-2& 1& 0& 0& 1 \\
\chi_{4c} & 4& 0& 2& 0& 1&-2& 0&-1& 0 \\
\chi_{4d} & 4& 0&-2& 0& 1&-2& 0& 1& 0 \\
\hline
\end{array} $$

\end{table}

\AbsT{Embedding $G'$ into $G$.}\label{subgrps}
Before proceeding to the local subgroups of $G$
we determine an explicit embedding of $G'$, and thus
of the whole subgroup chain $P<H'<G'$, into $G$.
Note that by \cite[Proof~of Corollary 6.5.5]{KleWil}
there is a unique conjugacy class of subgroups of $G$
isomorphic to $\fA_8$, and recalling that $G$ has
a unique conjugacy class of subgroups isomorphic to $C_3\tm C_3$ 
we may safely choose $P<G'$ as our favourite defect group of $A$.
In practice, we start with standard generators of $G$,
in the sense of \cite{Wil}, which hence in particular 
are unique up to simultaneous $G$-conjugacy.
Then we obtain straight-line programs 
producing generators of all of the subgroups in question.
We describe the computations which have to be done, 
using \GAP{} and the \MA{}:

\abs
In order to find a subgroup isomorphic to $\fA_8$ of $G$,
we work through the maximal subgroup $M:=2^{11}\cn M_{24}$ of $G$, 
since it also follows from \cite[Proof~of Corollary 6.5.5]{KleWil} 
that we even 
have a subgroup chain $G'<M_{24}<M<G$. Hence we 
start with the absolutely irreducible $\F_2$-representation $V$ of $G$ 
of dimension $112$, standard generators of which are available in 
\cite{ModAtlasRep}, together with a straight-line program yielding 
generators of $M$.
By \cite[p.190]{Atlas}, the restriction 
$\Res_M(V)$ is uniserial 
with radical series $[1a/11b/44b/44a/11a/1a]$.
Using the \MA{}, we compute the action of $M$ on the unique 
submodule $U$ of dimension $12$, which hence has shape $[11a/1a]$.
Considering the action of $M$ on the elements of $U$
yields a transitive permutation representation $\pi$ on $1518$ points, 
from which \GAP{} shows that we indeed have found a faithful
representation of $M$. 
Further computations in $M$ can now be done using the small 
representation $\pi$, where in particular we have all the 
facilities of \GAP{} to deal with permutation groups at our disposal.

\abs
Now we proceed similar to the approach in \cite[(16.2)]{MueHabil}:
Since the normal $2$-subgroup $2^{11}\lhd M$ acts trivially
on the constituent $11a$ of 
$\Res_M(V)$, 
we thus obtain a
representation of $M/2^{11}\cong M_{24}$, although in terms
of non-standard generators. Following the recipe given in \cite{ModAtlasRep},
by a random search we find a straight-line program producing
standard generators of the action of $M_{24}$ on $11a$, and
running this on the representation $\pi$, \GAP{} shows that 
we have indeed found standard generators of a subgroup $M_{24}$ of $M$.
Having this in place, we apply a straight-line program 
yielding generators of $2^4\cn\fA_8$ from standard generators 
of $M_{24}$, again available in \cite{ModAtlasRep}, and finally
another random search yields a straight-line program producing
generators of a subgroup $G':=\fA_8<2^4\cn\fA_8$.

\AbsT{Local subgroups of $G$.} 
Let $N:=N_G(P)$. 
It is stated in \cite[Section 3 and Proof~of Corollary 6.4.4]{KleWil},
unfortunately without proof, that $N\cong(2^3\tm P)\cn\GL_2(3)$ 
and that $N$ can be embedded into the maximal subgroup $M$. 
Since we need an explicit realisation of $N$ as a subgroup of $G$ anyway, 
we are going to verify these facts independently:

\abs
We have already fixed $P<G'<M$. We remark that,
since, by \cite[Proof~of Corollary 6.5.5]{KleWil}, 
$M$ has a unique conjugacy class of subgroups 
isomorphic to $\fA_8$, and $P<G'$ is a Sylow $3$-subgroup,
this defines a unique conjugacy class of subgroups of $M$ 
isomorphic to $C_3\tm C_3$.
Using \GAP{}, the normaliser $N_M(P)$ can be determined
explicitly as a subgroup of $M$.
It turns out that $|N_M(P)|=3456=2^3\cdot 3^2\cdot 48$.
We already know that $C=C_G(P)\cong 2^3\tm P$, hence from
$\Out(P)\cong\Aut(P)\cong\GL_2(3)$ we conclude that $C<M$
and that $P$ is fully automised in $M$, thus we infer $N=N_M(P)$.
Moreover, it is verified using \GAP{} that $C$ has a complement in $N$.

\abs
Let $B$ be the block of $kN$ which is the Brauer 
correspondent of $A$. Since $A$ is the only block of $kG$ 
having $P$ as a defect group, by Brauer's First Main Theorem, 
see \cite[Theorem 5.2.15]{NT}, we conclude that $B$
is the only block of $kN$ having $P$ as a defect group.

\AbsT{Local subgroups of $G$, continued.}
Next we proceed to find the blocks of $C$ covered by $B$:
Let $2^3\cong E<C$ be the elementary abelian Sylow $2$-subgroup.
Then then blocks of $C$ are in bijection with the irreducible
ordinary characters of $E$. 
Identifying $E$ with the $\F_2$-vector space $\F_2^3$, 
and its character group with 
$(\F_2^3)^\ast:=\Hom_{\F_2}(\F_2^3,\F_2)$,
the \MA{} shows that $(\F_2^3)^\ast$, with respect to the action
of $N/C\cong\GL_2(3)$, becomes a uniserial module 
with radical series $[2/1]$, where $\GL_2(3)$ has three orbits,
of lengths $[1,1,6]$, on the elements of $(\F_2^3)^\ast$.

\abs
Hence there are two irreducible characters of $E$ which are $N$-invariant, 
thus extend to linear
characters of $N$, and hence belong to blocks of (maximal) defect $3$ of $N$. 
Thus we conclude that $B$ covers an $N$-orbit of six blocks of $C$,
hence the associated inertia groups have index $6$ in $N$. 
Now we make a specific choice amongst these blocks:
 
\abs
We observe that $H'C/C\cong H'/(H'\cap C)=H'/C'\cong D_8$,
hence $H'C$ has index $6$ in $N$. This is sufficiently suspicious
to be tempted to check the following: Indeed, 
the \MA{} shows that $(\F_2^3)^\ast$, with respect to the action
of $H'C/C\cong D_8$, has radical series $[1/(1\oplus 1)]$, 
where $D_8<\GL_2(3)$ has five orbits, of lengths $[1,1,1,1,4]$, 
on the elements of $(\F_2^3)^\ast$. 
This says that there are precisely two blocks
of $kC$ amongst those covered by $B$ which are fixed by $H'C$. 
Moreover, using the identification of the elements of 
$(\F_2^3)^\ast$ with the irreducible ordinary characters of $E$
yields the associated primitive idempotents 
$e_i\in kE$, for $i\in\{1,2\}$, which are the block idempotents in $kC$
of the blocks in question. 

\abs
Letting $H_i:=N_G(P,e_i):=\{g\in N_G(P);e_i^g=e_i\}$ 
be the associated inertia groups, this shows that we have
$H:=H'C=H_1=H_2\cong (2^3\tm P)\cn D_8$.
This completes the panorama of subgroups of $G$ we need,
the complete picture is shown in Table \ref{J4subgrptbl}.
Finally, let $B_i$ be the block of $kH$ being the 
Fong-Reynolds correspondent of $B$ with respect to $e_i(kC)e_i$, 
that is $B_i$ covers $e_i(kC)e_i$ and we have $(B_i)^N=B$,
see \cite[Theorem 5.5.10]{NT}; note that we hence have $1_{B_i}=e_i$.

\begin{table}\caption{Subgroups of $G$.}\label{J4subgrptbl}
\begin{center} \unitlength=1pt
\boxed{\begin{picture}(338,366)(-100,-327)
\put(40,25){\framebox(60,12){$G=J_4$}}
\put(10,-25){\framebox(120,12){$M=2^{11}\cn M_{24}$}}
\put(-100,-100){\framebox(160,12){$N=N_G(P)=(2^3\tm P)\cn\GL_2(3)$}}
\put(178,-125){\framebox(60,12){$G'=\fA_8$}}
\put(-100,-150){\framebox(160,12){$H=N_G(P,e)=(2^3\tm P)\cn D_8$}} 
\put(100,-200){\framebox(120,12){$H'=N_{G'}(P)=P\cn D_8$}}
\put(-80,-225){\framebox(120,12){$C=C_G(P)=2^3\tm P$}}
\put(10,-275){\framebox(120,12){$C_{G'}(P)=P=C_3\tm C_3$}}
\put(40,-325){\framebox(60,12){$\{1\}$}}
\put(70,25){\line(0,-1){38}}
\put(70,-25){\line(-4,-3){83}}
\put(70,-25){\line(3,-2){131}}
\put(-20,-100){\line(0,-1){38}}
\put(-20,-150){\line(0,-1){63}}
\put(47,-150){\line(3,-1){113}}
\put(-20,-225){\line(3,-2){56}}
\put(208,-125){\line(-3,-4){47}}
\put(152,-200){\line(-3,-4){47}}
\put(70,-275){\line(0,-1){38}}
\put(-25,-119){\makebox(0,0){$6$}}
\put(-25,-182){\makebox(0,0){$8$}}
\put(102,-162){\makebox(0,0){$8$}}
\put(0,-247){\makebox(0,0){$8$}}
\put(136,-232){\makebox(0,0){$8$}}
\put(75,-295){\makebox(0,0){$9$}}
\put(197,-157){\makebox(0,0){$280$}}
\put(163,-66){\makebox(0,0){$24\,870\,912$}}
\put(42,7){\makebox(0,0){$173\,067\,389$}}
\put(-2,-54){\makebox(0,0){$145\,080\,320$}}
\end{picture}} \end{center}

\end{table}

\AbsT{Local structure of $G$.}\label{localstruct}
To determine the decomposition matrices of $B_i$ and $B$ 
we proceed as follows:
Since $H$ is a split extension of $E\cong 2^3$ by $H'\cong P\cn D_8$, 
the irreducible ordinary and Brauer characters of $H$ can be
determined from those of $E$ and $H'$ via Clifford theory, 
see \cite[Chapter 3.3]{NT}. Moreover, since $B_i$ covers $e_i(kE)e_i$, 
by \cite[Lemma 5.5.7]{NT} we are only 
interested in the characters of $H$ covering the 
irreducible character $\lb_i$ of $E$ associated with $e_i$. 
Since $\lb_i$ is $H$-invariant, it extends to $H$
by letting $\lb_i(gh):=\lb_i(g)$, for all $g\in E$ and $h\in H'$.
Then the irreducible ordinary and Brauer characters of $B_i$ 
are in natural bijection to those of $kH'=B'$.
Hence, using the notation in \ref{hprime}, 
we may write the irreducible ordinary characters of $B_i$
as $\chi_{1x_i}$, $\chi_{2_i}$, and $\chi_{4x_i}$,
and the simple $B_i$-modules as $1x_i$ and $2_i$, 
where $x\in\{a,b,c,d\}$, subject to the conditions 
$\Res_{H'}^H(\chi_{?_i})=\chi_?$ and $\Res_{H'}^H(?_i)\cong ?$
for characters and modules, respectively.
Thus, using these identifications, the decomposition matrices of 
$B_i$ and $B'$ coincide, see Table \ref{a8dectbl}.

%

\abs
Moreover, by Fong-Reynolds correspondence, see \cite[Theorem 5.5.10]{NT},
we may write the irreducible ordinary characters of $B$
as $\chi_{6x}$, $\chi_{12}$, and $\chi_{24x}$,
and the simple $B$-modules as $6x$ and $12$, 
where $x\in\{a,b,c,d\}$, subject to the conditions
$\Ind_H^N(\chi_{?_i})=\chi_{6\cdot ?}$ and $\Ind_H^N(?_i)\cong 6\cdot ?$.
Since $\Res_{H'}^H(\Res_H^N(\chi_{6\cdot ?})\cdot 1_{B_i})=\chi_?$
and $\Res_{H'}^H(\Res_H^N(6\cdot ?)\cdot 1_{B_i})\cong ?$,
the parametrisation of characters and modules of $B$
is indeed independent from the choice of Fong-Reynolds correspondent
$B_i$ we work through. Thus, using these identifications, the 
decomposition matrices of $B$ and $B'$ coincide as well, 
see Table \ref{a8dectbl}.

\section{The final stroke}\label{stroke}

\AbsT{Equivalence between $A'$ and $B'$.}\label{equivaprimebprime}
By \cite[Lemma 2.4]{KosKunWak}, there is a unique indecomposable 
direct summand $\cM'$ of the $kG'$-$kH'$-bimodule 
$1_{A'}(kG')1_{B'}$ with vertex 
$\Dt P:=\{(g,g)\in G'\tm H';g\in P\}<G'\tm H'$. 
Then, letting $\modcat{}{A'}$ and $\modcat{}{B'}$ denote 
the categories of finitely generated right $A'$- and $B'$-modules, 
respectively, by \cite[Lemma 6.5(iii)]{KosKunWak} we have a splendid 
stable equivalence of Morita type given by the tensor functor 
$?\otm_{kG'}\cM'\cn\modcat{}{A'}\ra\modcat{}{B'}$;
to shorten notation we just write $\cM'$ instead.
In particular, by \cite[Theorem 2.1(ii)]{Lin}, the functor $\cM'$ maps
any simple $A'$-module to an indecomposable $B'$-module.

\abs
Let $f'$ be the Green correspondence
with respect to the triple $(G',P,H')$.
Then, by \cite[Lemma A.3]{KosMueNoe},
for any indecomposable $A'$-module $V$ having $P$ as a vertex, 
the unique non-projective indecomposable summand of $\cM'(V)$ 
coincides with the Green correspondent $f'(V)$. Thus, recalling that, 
by \cite[Corollary 3.7]{Kno}, any simple $A'$-module $T$ indeed has 
$P$ as a vertex, we conclude that $\cM'(T)\cong f'(T)$.

\abs
The Green correspondents of the
simple $A'$-modules are known by \cite[Example 4.3]{Oku},
and given at the right hand side of Table \ref{functbl}.
Note that, using the notation in \ref{a8} and \ref{hprime},
we have $f'(1)\cong 1a$ and $f'(7)\cong 1b$, while $\{1c,1d\}$ 
are now distinguished by fixing the embedding $H'<G'$ and specifying $f'$.
Nowadays it is easy to verify this independently, 
by computing the above Green correspondents explicitly,
using \GAP{} and the \MA{}. 
Moreover, in the spirit of 
the present paper, here we already encounter a nice toy application 
of the Brauer construction, which we cannot resist to include:

\abs
We consider the natural permutation action
of $G'$ on $\Om:=\{1\ld 8\}$, hence
$k[\Om]\cong 1\oplus 7$ as $k[G']$-modules.
We may assume that
$P:=\lr{(1,2,3),(4,5,6)}<\Stab_{G'}(8)\cong\fA_7$, hence $\Om^P=\{7,8\}$. 
Since $P<\fA_7<G'$ is a Sylow $3$-subgroup,
by \ref{brauertrivsrc} and \ref{brauerperm} we get 
$k[\Om](P)\cong k[\Om^P]\cong k[N_{\fA_7}(P)\bsl H']$ 
as $k[H']$-modules. Hence from $P\cn 4\cong N_{\fA_7}(P)\lhd H'$ we infer 
$$ k[\Om](P)\cong k[H'/(P\cn 4)]\cong k[D_8/C_4]\cong 1a\oplus 1b .$$
Thus, it follows from \ref{puigthm} that
$f'(1)\oplus f'(7)\cong 1a\oplus 1b$ as $k[H']$-modules,
and since $f'(1)\cong 1a$ anyway, we get $f'(7)\cong 1b$.

%

\AbsT{Equivalence between $A$ and $B$.}
To relate the blocks $A$ and $B$ we work through the
Fong-Reynolds correspondents $B_i$, taking both of them into account:
By \cite[Lemma 2.4]{KosKunWak}, there is a unique indecomposable
direct summand $\cM_i$ of the $kG$-$kH$-bimodule
$1_A(kG)1_{B_i}$ with vertex $\Dt P<G\tm H$,
which by \cite[Lemma 6.3(i)]{KosKunWak}
induces a splendid stable equivalence of Morita type 
$?\otm_{kG}\cM_i\cn\modcat{}{A}\ra\modcat{}{B_i}$,
which we abbreviate by writing $\cM_i$.

\abs 
Moreover, arguing as in \cite[Theorem 1.5]{KosKunWak2004},
it follows from Fong-Reynolds correspondence, see
\cite[Theorem 5.5.10]{NT}, and \cite[Theorem 0.2]{Bro1990} that
the $kN$-$kH$-bimodule $\cN_i:=1_B(kN)1_{B_i}$
and its dual, the $kH$-$kN$-bimodule $\cN_i^\vee:=1_{B_i}(kN)1_B$,
induce a pair of mutually inverse Puig equivalences given by 
$?\otm_{kN}\cN_i\cn\modcat{}{B}\ra\modcat{}{B_i}\cn
 V\mt\Res_H^N(V)\cdot 1_{B_i}$
and
$?\otm_{kH}\cN_i^\vee\cn\modcat{}{B_i}\ra\modcat{}{B}\cn V\mt\Ind_H^N(V)$,
which we abbreviate by $\cN_i$ and $\cN_i^\vee$, respectively.
Hence by concatenation we get a splendid stable equivalence of Morita type
$\cN_i^\vee\circ\cM_i\cn\modcat{}{A}\ra\modcat{}{B}$,
again mapping any simple $A$-module to an indecomposable $B$-module. 

\abs
We are now able to relate $\cN_i^\vee\circ\cM_i$ to
the Green correspondence $f$ with respect to the triple $(G,P,N)$:
The functor $\cN_i^\vee\circ\cM_i$ coincides with the tensor functor
$?\otm_{kG}(\cM_i\otm_{kH}\cN_i^\vee)\cn\modcat{}{A}\ra\modcat{}{B}$,
where the 
$kG$-$kN$-bimodule $\cM_i\otm_{kH}\cN_i^\vee$ is a direct summand of
$1_A(kG\otm_{kH}kN)1_B\cong 
 1_A\left(\bigoplus_{i=1}^{[N\cn H]}kG\right)1_B$.
Hence from \cite[Lemma A.3]{KosMueNoe} we infer that,
for any indecomposable $A$-module $V$ having $P$ as a vertex,
the unique non-projective indecomposable summand of $\cN_i^\vee(\cM_i(V))$
coincides with the Green correspondent $f(V)$;
see \cite[Lemma 6.3(iii)]{KosKunWak}.

\abs
Thus, recalling that, by \cite[Corollary 3.7]{Kno},
any simple $A$-module $S$ indeed has $P$ as a vertex,
we conclude that $\cN_i^\vee(\cM_i(S))\cong f(S)$,
in particular saying that $\cN_1^\vee\circ\cM_1$
and $\cN_2^\vee\circ\cM_2$ coincide on simple $A$-modules.
Hence $(\cN_2^\vee\circ\cM_2)^{-1}\circ(\cN_1^\vee\circ\cM_1)$
is a stable auto-equivalence of $\modcat{}{A}$, which is
the identity on simple $A$-modules, thus by \cite[Theorem 2.1(iii)]{Lin} 
is equivalent to the identity functor on $\modcat{}{A}$.
Hence we have 
$\cM:=\cN_1^\vee\circ\cM_1\cong\cN_2^\vee\circ\cM_2\cn
 \modcat{}{A}\ra\modcat{}{B}$.

\abs
From the equivalence $\cM_i\cong\cN_i\circ\cM$  
we get $\cM_i(S)\cong\Res_H^N(f(S))\cdot 1_{B_i}$,
for any simple $A$-module $S$.
Moreover, by \cite[Lemmas 6.7--6.10]{KosKunWak}
the modules $\cM(S)$ and $\cM_i(S)$ are as given in Table \ref{functbl},
where $\{\al,\bt,\gm,\dt\}=\{a,b,c,d\}$.
Hence it remains to connect the left and right hand sides of 
Table \ref{functbl}:

\AbsT{Equivalence between $B_i$ and $B'$.}\label{equivbbprime}
Similar to the equivalence between $B$ and $B_i$ above, 
we now consider the $kH$-$kH'$-bimodule $\cL_i:=1_{B_i}(kH)1_{B'}$
and its dual, the $kH'$-$kH$-bimodule $\cL_i^\vee:=1_{B'}(kH)1_{B}$.
Then it follows from \ref{localstruct} and \cite[Theorem 0.2]{Bro1990} that
$?\otm_{kH}\cL_i\cn\modcat{}{B_i}\ra\modcat{}{B'}\cn V\mt\Res_{H'}^H(V)$,
and
$?\otm_{kH'}\cL_i^\vee\cn\modcat{}{B'}\ra\modcat{}{B_i}\cn
 V\mt\Ind_{H'}^H(V)\cdot 1_{B_i}$ 
are a pair of mutually inverse Puig equivalences,
which we abbreviate by $\cL_i$ and $\cL_i^\vee$, respectively.

\abs
By \cite[Lemma 2.8]{KosKunWak}, tensoring with a linear modules
induces a Puig auto-equivalences
$?\otm_k 1z\cn\modcat{}{B'}\ra\modcat{}{B'}\cn V\mt V\otm 1z$,
for $z\in\{a,b,c,d\}$. 
Here, the trivial module $1a$ induces the identity functor on $\modcat{}{B'}$, 
while, the group of linear characters of $H'\cong P\cn D_8$
being isomorphic to $H'/[H',H']\cong C_2\tm C_2$,
for $z\neq a$ we get non-trivial involutory auto-equivalences.
This yields Puig equivalences
$\cL_i^z\cn\modcat{}{B_i}\ra\modcat{}{B'}\cn V\mt\Res_{H'}^H(V)\otm_k 1z$,
where of course $\cL_i\cong\cL_i^a$; in particular we get Puig equivalences
$\cL_i^z\circ\cN_i\cn\modcat{}{B}\ra\modcat{}{B'}$.

\abs
Moreover, twisting with the non-inner automorphism $\om\in\Aut(H')$, 
see \ref{hprime}, induces an involutory Morita auto-equivalence 
$\cW\cn\modcat{}{B'}\ra\modcat{}{B'}\cn V\mt V^\om$. Since 
applying $\om$ changes the embedding of $P$ into $H'$, this
the functor $\cW$ is \emph{not} a Puig auto-equivalence; 
see \cite[Lemma 6.12]{KosKunWak}.

\abs
Hence, to complete the picture in Table \ref{functbl}, we apply 
$\cL_i^\al$, where $1\al:=\Res_{H'}^H(1\al_i)$;
note that despite notation $\al\in\{a,b,c,d\}$ depends on $i\in\{1,2\}$.
Then we get $\cL_i^\al(1\al_i)=1a$, and hence
it follows from \cite[Lemmas 6.8 and 6.10]{KosKunWak} that $\cL_i^\al(1\bt_i)=1b$.
This implies $\cL_i^\al(\{1\gm_i,1\dt_i\})=\{1c,1d\}$, where 
$1d\cong (1c)^\om$.
Thus we get a stable equivalence of Morita type
$$\cF_i:=\cM^{\prime -1}\circ\cW^\eps\circ\cL_i^\al\circ\cN_i\circ\cM
  \cn\modcat{}{A}\ra\modcat{}{A'} ,$$ 
which for an appropriate choice of $\eps\in\{0,1\}$,
again depending on $i\in\{1,2\}$, maps simple $A$-modules 
to simple $A'$-modules, hence by \cite[Theorem 2.1(iii)]{Lin} 
is an equivalence; note that, by \cite[Lemma 4.6.(iii)]{KosKunWak},
the uniserial modules appearing in Table \ref{functbl} 
are uniquely determined up to isomorphism by their radical series.

\abs
Moreover, by construction, $\cF_i$ is a splendid equivalence,
that is a Puig equivalence, if and only if $\eps=0$. Thus, in
our setting, the question left open in \cite[Question 6.14]{KosKunWak}
can be reformulated as follows: Is $\eps=0$ or $\eps=1$?
Recall that we still have two cases $i\in\{1,2\}$ at our disposal.
Since there does not seem to be a way to answer this by abstract 
theory alone, we force a decision by explicit computation,
which finally is our main application of the Brauer construction:

\begin{table}\caption{Morita equivalence between $A$ and $A'$}\label{functbl}
$$ \begin{array}{|ccccccccc|}
\hline
\modcat{}{A} & \overset{\cM}{\lra} & \modcat{}{B}
             & \overset{\cN_i}{\lra} & \modcat{}{B_i}
             & \overset{\cL_i^\al}{\lra} & \modcat{}{B'}
             & \overset{\cM'}{\lla} & \modcat{}{A'} \\
\hline\hline
S_1 & \mt & 6\al & \mt & 1\al_i & \mt & 1a
    & \mapsfrom & 1 \rule{0em}{1.5em} \\
S_3 & \mt & 6\bt & \mt & 1\bt_i & \mt & 1b 
    & \mapsfrom & 7 \rule{0em}{1.5em} \\
S_5 & \mt & \boxed{\begin{matrix}6\dt \\ 12 \\ 6\dt\end{matrix}}
    & \mt & \boxed{\begin{matrix}1\dt_i \\ 2 \\ 1\dt_i\end{matrix}}
    & \overset{\cW^\eps}{\mt} 
    & \boxed{\begin{matrix}1d \\ 2 \\ 1d\end{matrix}}
    & \mapsfrom & 13 \rule{0em}{3em} \\
S_2 & \mt & 6\gm & \mt & 1\gm_i & \overset{\cW^\eps}{\mt} & 1c
    & \mapsfrom & 28 \rule{0em}{1.5em} \\
S_4 & \mt & 12 & \mt & 2_i & \mt & 2 
    & \mapsfrom & 35 \rule{0em}{1.5em} \\
\hline
\end{array} $$

\end{table}

\AbsT{Applying the Brauer construction.}\label{applybrauer}
We proceed to determine the Green correspondents 
$f(S_1)\cong\cM(S_1)\cong 6\al$ and $f(S_2)\cong\cM(S_2)\cong 6\gm$
of $S_1$ and $S_2$. Recall that $S_1$ and $S_2$ are of dimension
$4\,290\,927$ and $95\,288\,172$, respectively, hence we cannot
possibly simply proceed as in \ref{equivaprimebprime}. Instead,
we consider the $G$-action on the set $\Om:=M\bsl G$ of cosets of 
$M$ in $G$, which has cardinality $173\,067\,389$:

\abs
Let $1_M^G := \Ind_M^G(1_M)$ denote 
the permutation character of the $G$-action on $\Om$.
Using the ordinary character tables of $M$ and $G$, contained
in the character table library of \GAP{}, it turns out that the
component of $1_M^G$ belonging to the block $A$ is given as
$(1_M^G)\cdot 1_A=\chi_{14}+\chi_{21}$; see \cite[Lemma 3.12]{KosKunWak}.
Thus, using the decomposition matrix of $A$, see \ref{j4}, we get
$k[\Om]\cdot 1_A\cong S_1\oplus S_2$ as $kG$-modules, hence both 
$S_1$ and $S_2$ are trivial source modules; see \cite[Lemma 3.13]{KosKunWak}.
Moreover, by \ref{puigthm} we have
$$ f(S_1)\oplus f(S_2)\cong S_1(P)\oplus S_2(P)\spmid k[\Om](P) .$$
Hence, by \ref{brauertrivsrc}
we proceed to determine $k[\Om]^P$ and its structure as a
permutation $kN$-module.
In order to apply \ref{brauerperm}, we have to find a set of
representatives of the $M$-conjugacy classes of subgroups
of $M$ being $G$-conjugate to $P$:

\abs
Using \GAP{}, a Sylow $3$-subgroup of $M$, being isomorphic to an
extraspecial group $3_+^{1+2}$,
and from that the conjugacy classes of subgroups of $M$ isomorphic
to $C_3\tm C_3$, can be determined. It turns out that $M$ has
precisely two such conjugacy classes. One of them of course containing $P$,
let $\wti{P}<M$ be a representative of the other conjugacy class.
Noting that $k[N_M(P)\bsl N_G(P)]=k[N\bsl N]\cong k$ 
is the trivial $kN$-module, from \ref{brauerperm} we thus get,
as permutation $kN$-modules,
$$ k[\Om^P]\cong k\oplus k[N_{M^g}(P)\bsl N] 
           \cong k\oplus k[N_M(\wti{P})\bsl N_G(\wti{P})]^g ,$$
for some $g\in G -  M$ such that $\wti{P}^g=P$.
Note that within our setting, only allowing to compute efficiently in $M$,
we are not able to get hands on a conjugating element easily;
hence we circumvent an explicit choice of such an element:

\abs
To determine the action of $N$ on the set $N_{M^g}(P)\bsl N$, 
up to equivalence of permutation actions,
it suffices to find a subgroup
of $N$ which is $N$-conjugate to $N_{M^g}(P)$. To this end,
employing \GAP{} again, we first compute $N_M(\wti{P})$ 
as a subgroup of $M$.
It turns out that we have $N_M(\wti{P})\cong 2\tm(P\cn D_{12})$,
where $D_{12}\cong 3\cn 2^2$ can be identified
with a Borel subgroup of $\GL_2(3)\cong\Aut(P)$, and using the 
group library of \GAP{} we get
$N_M(\wti{P})\cong\textsf{SmallGroup}(216,102)$.
Now we have $N_M(\wti{P})\cong N_{M^g}(P)<N$, and
\GAP{} shows that $N$ has a unique conjugacy class
of subgroups isomorphic to $\textsf{SmallGroup}(216,102)$.
Hence letting $\wti{N}<N$ be a representative of this conjugacy class,
we have $k[N_{M^g}(P)\bsl N]\cong k[\wti{N}\bsl N]$ 
as permutation $kN$-modules, and thus 
$k[\Om^P]\cong k\oplus k[\wti{N}\bsl N]$.

\abs
Note that, as far as we see, we are just lucky here:
If there were several conjugacy classes of subgroups of $N$
of the above isomorphism type, then it might become necessary
to construct a conjugating element in 
$G - M$ explicitly.
Anyway, although we have started with the huge $G$-permutation domain $\Om$,
we have now managed to reduce the problem to a consideration
of the tiny $N$-permutation domain $\wti{N}\bsl N$; 
note that $[N\cn\wti{N}]=16$:

\AbsT{Conclusion.}
Using \GAP{}, we determine the permutation action of $N$ on the cosets
of $\wti{N}$. Then it is straightforward, using the \MA{}, to find the
constituents of the permutation module $k[\wti{N}\bsl N]$.
It turns out that there are precisely two constituents 
of dimension $6$,
$6_I$ and $6_{II}$ say, each occurring 
with multiplicity one. 
Hence by \ref{brauertrivsrc} and \ref{brauerperm} we have 
$\{f(S_1),f(S_2)\}=\{6_I,6_{II}\}$.
Applying $\cL_1\circ\cN_1$ and $\cL_2\circ\cN_2$ to $6_I$ and $6_{II}$, 
the \MA{} shows that
$$ \begin{array}{lcl}
\Res_{H'}^H(\Res_H^N(6_I)   \cdot(1_{B_1}+1_{B_2}))& \cong & 1a\oplus 1a,\\ 
\Res_{H'}^H(\Res_H^N(6_{II})\cdot(1_{B_1}+1_{B_2}))& \cong & 1c\oplus 1d.
                                                     \rule{0em}{1.5em} \\
\end{array} $$
Hence we may assume that
$\Res_{H'}^H(\Res_H^N(6_{II})\cdot 1_{B_1})\cong 1c$.

\abs
Thus, using the notation from \ref{equivbbprime},
the functor $\cF_1$ is equipped with the
parameter $\al=a$ or $\al=c$, depending on whether
$f(S_1)\cong 6_I$ or $f(S_1)\cong 6_{II}$, 
where in both cases we have $\eps=0$,
thus $\cF_1$ is a Puig equivalence.
With respect to the same case distinction, $\cF_2$ is equipped with the
parameter $\al=a$ or $\al=d$, 
where in both cases we have $\eps=1$,
hence $\cF_2$ is \emph{not} a Puig equivalence.
Anyway, this proves Theorem  \ref{puigequiv}.

\Abs{\bf Conjecture/Question.}
Finally, we conjecture that, due to our fixed subgroup
configuration, we have $f(S_1)\cong 6_I$, that is
$\Res_{H'}^H(\Res_H^N(f(S_1))\cdot 1_{B_i})\cong 1a$ for both $i\in\{1,2\}$.
Again this cannot be answered by abstract theory alone:
Indeed, \GAP{} shows that $H$ has four conjugacy classes
of subgroups isomorphic to $P\cn D_8$, all of which supplement $C$,
hence each can be used to parametrise the characters of $B_i$
along the lines of \ref{localstruct}, but the resulting parametrisations 
are distinct. Unfortunately, our explicit computations, as far as we have
pursued them, do not tell us either.


\medskip\noindent
{\bf Acknowledgements.}
\small{\rm
This work was done while the first author was
staying in RWTH Aachen University in 2011 and 2012. 
He is grateful to Gerhard Hiss for his kind hospitality.
For this research the first author was partially
supported by the Japan Society for Promotion of Science (JSPS),
Grant-in-Aid for Scientific Research (C)23540007, 2011--2014. 
The second author is grateful for financial support in
the framework of the DFG (German Science Foundation) Priority
Programme SPP-1388 `Representation Theory', which this
research is a contribution to.}



\end{document}